# Comprehensive Review of Analytical and Numerical Approaches in Earth-to-Air Heat Exchangers and Exergoeconomic Evaluations

Saeed Asadi[1, *], Mohsen Mohammadagha[1], Hajar Kazemi Naeini[1]

[1]Department of Civil Engineering, University of Texas at Arlington, Arlington, Texas.

[*]Corresponding author: sxa1930@mavs.uta.edu

## Abstract

In recent decades, Earth-to-Air Heat Exchangers (EAHEs), also known as underground air ducts, have garnered significant attention for their ability to provide energy-efficient cooling and heating solutions while maintaining a minimal environmental footprint. These systems leverage the relatively stable underground temperature to regulate indoor climates, reducing reliance on conventional heating, ventilation, and air conditioning (HVAC) systems. This review systematically categorizes and synthesizes research on EAHEs into three primary areas: analytical, numerical, and exergoeconomic studies. Analytical approaches focus on developing theoretical models to predict thermal performance, while numerical simulations provide insights into system optimization and real-world applications. Exergoeconomic analyses, integrating thermodynamic efficiency with economic considerations, offer valuable perspectives on cost-effectiveness and long-term viability. By consolidating existing contributions across these domains, this study serves as a comprehensive reference for researchers, engineers, and policymakers seeking to enhance the design, implementation, and performance of EAHE systems. The findings emphasize the pivotal role of EAHEs in reducing energy consumption, lowering greenhouse gas emissions, and improving economic sustainability. Additionally, this review identifies key challenges, including soil thermal conductivity variations, moisture effects, and system integration with renewable energy sources, which require further investigation. By addressing these challenges, EAHEs can be further optimized to serve as a cornerstone in sustainable energy management, contributing to global efforts toward energy-efficient building solutions and climate change mitigation.

# 1- Introduction

In 1997, Kyoto, Japan, hosted the International Seminar on Climate Change. The seminar's conclusions emphasized the shared responsibility of all nations, whether developed or developing, to mitigate environmental damage caused by the improper use of non-renewable energy sources. Countries were tasked with adhering to principles that preserve environmental health. Climate change, global warming, ozone layer depletion, and acid rain are all consequences of inefficient energy use and excessive carbon dioxide emissions. This undeniable reality has prompted humanity to seek solutions to achieve the objectives outlined during the Kyoto seminar. The primary concern is not solely economic stability though it is crucial but rather the protection of the environment for both current and future generations. Unfortunately, energy consumption practices in developing countries, including Iran, often highlight a lack of attention to these critical issues. For instance, the growing use of air conditioning systems, particularly in office buildings, underscores the energy inefficiency in these regions.

These systems, designed to maintain thermal comfort in buildings, consume substantial amounts of energy, often during peak electricity demand periods. In regions with extreme summer heat, providing thermal comfort via air conditioning systems becomes exceedingly costly and often inefficient. For example, Camille (1988) found that 59% of electricity consumption in Riyadh was used for cooling buildings [1]. Similarly, Saeed (1997) reported that 37% of the annual electricity production in Riyadh was consumed for cooling during the hot seasons [2]. A scientific evaluation of the benefits and drawbacks of using such systems reveals that their limited advantages often fail to justify their widespread use. Modern technology has introduced solutions for achieving thermal comfort in buildings without considering climatic conditions. Today, in many countries, air conditioning systems are perceived as a measure of building quality, especially in commercial and office spaces. In the UK, for instance, 100% of new buildings in 1999 were equipped with air conditioning systems, despite maximum summer temperatures rarely exceeding 25°C.

Energy consumption in such buildings is at least four times higher than in others, leading to a proportional increase in carbon dioxide emissions, which pose significant environmental challenges. The 20% rise in air conditioning system sales in Greece, Spain, and Italy compared to 1999 is another indication of the disregard for environmental and economic sustainability. Additionally, maintenance costs for these systems over a 10-year period, adjusted for inflation, often equal the initial purchase and installation costs. Buildings equipped with air conditioning systems are also prone to humidity-related issues. The most significant concern, however, is the high electricity consumption of these systems during peak demand periods, often coinciding with the hottest months and hours of the year. This elevated demand strains natural and atmospheric resources and threatens energy reserves. Why do engineers and building experts fail to address these issues? A lack of consideration for national interests, environmental health, and future generations, coupled with the influence of global markets, has perpetuated the widespread use of these systems.

Admittedly, certain regions with unique climatic conditions may require air conditioning systems. However, it is essential to first identify these areas and evaluate whether such systems are necessary. If needed, efforts should focus on determining the extent of this necessity and exploring alternative solutions or energy-efficient systems to address the problem. Poorly designed buildings that fail to respond to climatic conditions will only increase demand for mechanical air conditioning systems. Kumar (1998) examined the impact of appropriate architectural designs on cooling and heating needs and thermal comfort [3]. His findings demonstrated that well-designed buildings tailored to local climates could reduce reliance on mechanical air conditioning systems. One method to minimize the demand for mechanical air conditioning systems is to employ natural ventilation systems. For instance, in Tunisia and eastern Spain, where the climate is hot and arid, underground houses were constructed. Similarly, in northern China, large underground homes were built to escape the severe cold of winter. In areas unsuitable for excavation, ancient architects utilized alternative techniques to ensure thermal comfort. Examples include the famous wind towers of Saudi Arabia, Pakistan, and Iran, which showcase innovative passive cooling techniques [4].

As previously mentioned, the concept of utilizing the ground as a thermal reservoir or underground pipes has roots in ancient practices. However, the implementation of Earth-to-Air Heat Exchangers (EAHEs) began in developed countries due to their higher thermal efficiency compared to purely convective systems [5]. This report categorizes the research in this field into three main sections analytical studies, numerical studies, and exergoeconomic analyses:, which combine economic and exergy evaluations. Before exploring these studies, it is essential to address the significance of using EAHE systems and their advantages over other air-to-ground systems. These systems are known for their low energy consumption and their ability to reduce air pollution significantly due to the absence of refrigerants. By using air as the working fluid, they eliminate the need for specialized fluids, further simplifying the design and reducing maintenance and operational costs. Additionally, their straightforward design minimizes the need for additional equipment in occupied spaces, making them economically viable in the long term.

Despite these advantages, EAHE systems have certain limitations. The initial installation costs are relatively high, primarily due to the excavation and burial of pipes or ducts. However, these costs are offset over time as the system requires minimal maintenance and operational expenses. Another drawback is the potential transmission of fan noise through the ducts into indoor spaces, which can be disruptive. In highly humid climates, condensation within the pipes may reduce system efficiency, although this issue can be mitigated by employing drainage pumps, albeit at the cost of increased energy consumption. Furthermore, in closed-loop systems, air quality may degrade over time, necessitating the use of air filters to maintain acceptable indoor air standards. Figure 1 illustrates the two primary configurations of EAHE systems, namely the open-loop and closed-loop designs, each with its specific operational characteristics and applications.

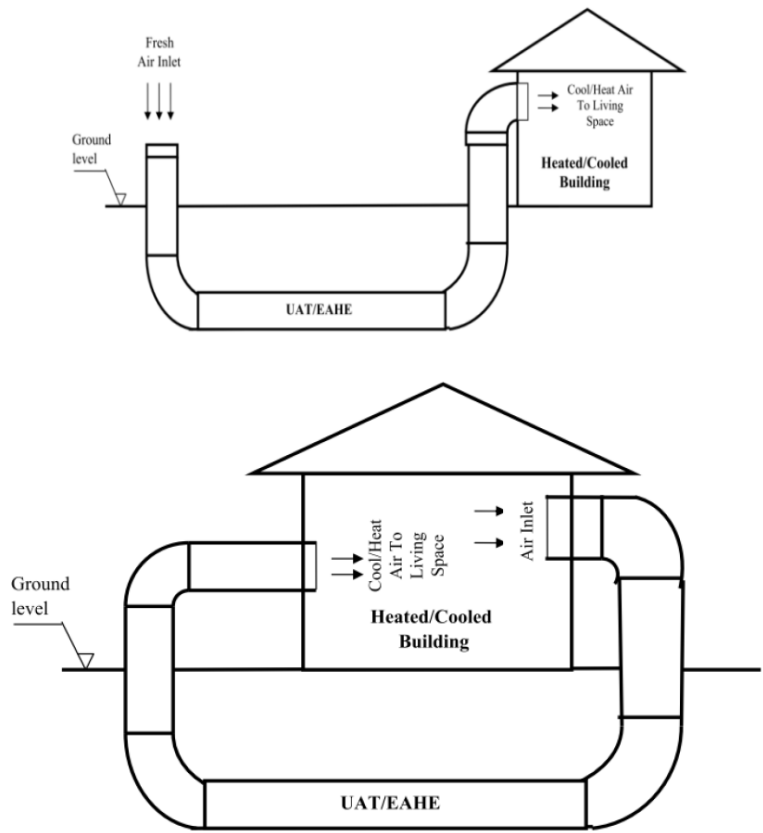

Figure 1: Schematic Design of Open and Closed Systems

## 2- Literature review

### *2-1- Analytical Studies Conducted*

Tzafriris and Liparakis (1992) reviewed eight algorithms proposed up to that time for predicting the performance of underground air ducts. Their study evaluated the sensitivity of each model to factors such as inlet air temperature, air velocity, duct length, diameter, and installation depth. By comparing experimental results with the models' outputs, they assessed the accuracy of each method [11]. Karate and Kreider (1995) developed a simple analytical model to determine the efficiency and thermal potential of underground air ducts. This model could predict the variations in the outlet air temperature over a 24-hour period. They used parametric analysis to examine the effects of hydraulic diameter and air flow rate on heat transfer. The results showed that increasing the hydraulic diameter had a more pronounced impact on improving thermal performance than increasing air velocity [12].

Jank et al. (1999) investigated heat transfer from buried pipes under the assumption of constant wall temperature using an analytical model. They employed the method of separation of variables to transform the problem from a circular domain in an infinite medium to a rectangular domain. By solving the governing differential equations numerically, they obtained the temperature distribution in the surrounding soil and the outlet air temperature [13]. Pap and Janssens (2003) analyzed the key parameters influencing the thermo-hydraulic design of underground ducts using a one-dimensional approximation and the NTU (Number of Transfer Units) method. Their analytical approach enabled the creation of design charts for underground ducts [14]. Kokomo et al. (2008) proposed a simplified one-dimensional analytical method to account for the condensation phenomenon in underground ducts. Their approach provided insights into the effects of condensation on system performance [15].

## 2-2-Numerical Studies

Al-Ajmi et al. (2006) used TRNSYS software to simulate a residential building in the hot and dry climate of Kuwait. They reported that while underground ducts alone could not provide full thermal comfort, their use for pre-cooling incoming air reduced the building's energy consumption during the cooling season by approximately 30% [16]. Zhang and Haghighat (2009) employed computational fluid dynamics (CFD) to study the thermal behavior of rectangular underground ducts. Their parametric analysis revealed that variations in surface temperature, inlet air turbulence, and outlet size had negligible effects on the Nusselt number [17]. Bansal et al. (2009) conducted a parametric study to evaluate the use of underground ducts for reducing heating loads in buildings during winter. A numerical model was used to assess the effects of parameters such as pipe material and air velocity on system performance. Their findings indicated that pipe material had minimal impact on thermal performance, suggesting that lower-cost materials could be used effectively [18].

Chau and Tiwari (2009) investigated the performance of a system combining an underground air duct and a dome-shaped roof for ventilating an indoor space. Their numerical and experimental analysis showed that this system reduced indoor air temperatures by 5–15°C in summer and increased them by a similar range in winter [19].

Joachim et al. (2011) examined the numerical and experimental use of an underground air duct for ventilating a 55-cubic-meter space. Using Fluent software and a three-dimensional simulation approach, they reported that the system cooled air by up to 4°C in summer and warmed it by up to 8°C in winter [20]. Asciune et al. (2011) studied the performance of underground air ducts in three different climatic zones in Italy. They found that for summer conditions, the best thermal performance was achieved with moist soil in the cold climate of Milan. They also concluded that duct lengths shorter than 10 meters or longer than 50 meters should be avoided and recommended a burial depth of approximately 3 meters [21]. Bansal et al. (2013) introduced a new parameter to describe the efficiency of underground air ducts more effectively. Using Fluent software, they studied the thermal performance of a system in India. Their results showed that the system's thermal efficiency declined over time after initial operation [22]. A decision-making framework incorporating Z-SWARA, Z-WASPAS, and FMEA was introduced by Lonbar et al. [28] to enhance budget allocation and risk evaluation in media advertising, thereby improving strategic planning quality. Asgari et al. [29] utilized threshold regression and Granger causality to analyze the connection between energy consumption and economic growth across ten OECD countries. Their research highlights that sustainable development efforts drive economic growth while ensuring environmental preservation.

Sharifi et al. [30] examined AI-based digital twin models for stormwater infrastructure in smart cities, emphasizing their significance in optimizing drainage systems. Their study suggests that standardized definitions and further research can enhance stormwater management and water quality. Sharifi et al. [31] applied gene expression programming to predict NOx emissions and mechanical efficiency in diesel engines, determining that the air/fuel ratio plays a crucial role. Correlation analysis and Sobol'-Jansen sensitivity assessments revealed nonlinear interactions influencing efficiency and emissions. Reihanifar et al. [32] explored nanotechnology's role in groundwater remediation, focusing on pollutant removal, sustainability challenges, and future developments. Their study underscores the importance of innovative nanomaterials and their environmental implications. Ahmadi et al. [33] developed a digital twin model integrating deep learning for coastal terrain segmentation. Utilizing U-Net and USGS data, their approach accurately classifies landscapes, aiding environmental monitoring and urban planning. Attari et al. [34] proposed an optimization model aimed at reducing logistics network delays and costs, employing meta-heuristic algorithms in MATLAB to efficiently distribute inventory across diverse storage conditions. Ahmadi et al. [35] conducted a comparative analysis of SAM and U-Net models for detecting cracks in concrete structures, demonstrating improved accuracy through SAM's strength in identifying longitudinal cracks and U-Net's precision in determining crack size and location. Sharifi et al. [36] used OpenFOAM to analyze geyser eruptions in storm sewer systems, investigating air-water interactions within an inverted Tee pipe and highlighting pressure's role in slug formation and geyser intensity. Moghim et al. [37] studied extreme hydrometeorological events in Bangladesh, assessing climate impacts on infrastructure reliability. Their findings improve rainfall simulations and enhance flood risk predictions. Asadi et al. [38] examined a passive cooling system that integrates a ground-to-air heat exchanger with a cylindrical roof, optimizing airflow and heat exchange efficiency through MATLAB and Fluent simulations. Their results demonstrate the system's effectiveness in regulating indoor temperatures.

Hussain et al. [39] employed conditional and random effects logistic regression models to analyze inter-seasonal meteorological drought dynamics in Ankara, Turkey. Their study highlights the influence of antecedent moisture conditions on drought transitions and underscores the utility of the RELogRM model for drought prediction. Cao et al. [40] introduced an optimized approach for airport security queue management using a Minimum Cost Network Flow (MCNF) model. Their solution enhances throughput, minimizes waiting times, and boosts operational efficiency, providing a scalable long-term strategy for airport security. Reihanifar et al. [41] designed a meteorological drought forecasting model that prioritizes both accuracy and simplicity, outperforming traditional forecasting methods based on Burdur City data. More et al. [42] presented the VMD-GP, an evolutionary model that improves meteorological drought prediction in ungauged catchments. Their method, which combines variational mode decomposition with genetic programming, enhances SPEI forecasting in Erbil, Iraq. Reihanifar et al. [43] applied value engineering to increase the sustainability of road construction projects, advocating for the integration of time, cost, and quality considerations to optimize project planning. Their research highlights value engineering's role in improving efficiency and reducing environmental impact. Reihanifar and Naimi [44] explored techniques to enhance quality and efficiency in project management, demonstrating how value management engineering can optimize material use, refine methodologies, and improve workplace relationships to boost productivity. Their study provides an overview of value management's principles, methodology, and global applications. Leon et al. [45] developed a finite volume model to enhance the stability of sewer simulations during filling and emptying phases. Their adaptation of the HLLS scheme ensures stationary conditions for both free surface and pressurized flows, thereby improving simulation accuracy. Mahyawansi et al. [46] conducted a large-scale experimental study on storm sewer geysers using an 88-meter-long pipeline setup. Their research examines the impact of the air-water volume ratio and pressure coefficient on geyser intensity.

## *2-3- Exergoeconomic Analyses of Earth-to-Air Heat Exchanger Systems*

This section reviews studies that have conducted exergy and economic analyses on systems equipped with Earth-to-Air Heat Exchangers (EAHEs). While EAHEs themselves do not consume energy directly, they often require auxiliary components such as fans to regulate airflow velocity, electric heaters for supplemental heating, and cooling devices to enhance performance under specific climatic conditions. The inclusion of these energy-consuming components necessitates exergy analysis to assess system efficiency and identify potential areas for optimization. Exergy analysis evaluates the quality and usability of energy within a system by determining the losses associated with irreversibilities, heat transfer inefficiencies, and pressure drops. By applying exergy-based methods, researchers can quantify the degradation of energy and pinpoint components with the highest exergy destruction. These insights are crucial for improving EAHE performance and integrating the system into broader energy-efficient frameworks.

In parallel, economic assessments of EAHE systems consider installation costs, operational expenses, and long-term financial benefits. Exergoeconomic analysis, which combines thermodynamic (exergy) principles with economic evaluation, provides a comprehensive approach to optimizing cost-effectiveness. Studies in this domain typically involve cost-to-benefit comparisons of EAHE-integrated HVAC systems, life cycle cost (LCC) analysis, and payback period assessments. The findings highlight that while EAHEs contribute to substantial energy savings over time, their economic viability depends on factors such as local climate conditions, soil thermal properties, and the efficiency of supporting equipment. Several studies have explored the potential of integrating EAHEs with renewable energy sources, such as solar-assisted heat pumps, to further enhance exergoeconomic performance. This hybrid approach reduces dependency on traditional HVAC systems, mitigates greenhouse gas emissions, and improves overall energy sustainability. However, challenges such as system complexity, initial investment costs, and maintenance requirements must be addressed to maximize the feasibility of widespread adoption. By synthesizing findings from existing literature, this section underscores the importance of exergoeconomic assessments in the advancement of EAHE technology. Future research should focus on refining optimization techniques, improving cost-effectiveness, and exploring innovative hybrid configurations to further enhance the sustainability and efficiency of EAHE systems.

## 3- Methodology

### 3-1- Mathematical Relations for Analysis

To provide a clearer understanding of the analytical approach, some of the mathematical relations utilized in these studies are briefly outlined [22–25].
General Exergy Balance Equation:

$$\dot{E}x_{in} - \dot{E}x_{out} = \dot{E}x_{dest} \tag{1}$$

Expanded Form of the Exergy Balance Equation:

$$\dot{E}x_{heat} - \dot{E}x_{work} + \dot{E}x_{mass,in} - \dot{E}x_{mass,out} = \dot{E}x_{dest} \tag{2}$$

If we expand the heat transfer at section *k* and the work rate in Equation (2), we arrive at Equation (3):

$$\sum \left(1 - \frac{T_0}{T_k}\right)\dot{Q}_k - \dot{W} + \sum \dot{m}_{in}\psi_{in} - \sum \dot{m}_{out}\psi_{out} = \dot{E}x_{dest} \tag{3}$$

Irreversibility or Exergy Destruction:

$$I = \dot{E}x_{dest} = T_0 \dot{S}_{gen} \tag{4}$$

$$\varepsilon^1 = \frac{\dot{E}x_{out}}{\dot{E}x_{in}} = 1 - \frac{\dot{E}x_{dest}}{\dot{E}x_{in}} \tag{5}$$

$$\varepsilon_{2,R} = \frac{\dot{E}x_{desired,output}}{\dot{E}x_{used}} \tag{6}$$

$$\varepsilon^3 = \frac{COP_{C,h}}{COP_{act,h}} \tag{7}$$

Where we have:

$$COP_{C,h} = \frac{T_H}{T_H - T_L} \tag{8}$$

Equation (8) represents the maximum coefficient of performance for an ideal system operating between a low temperature T*l* and a high temperature T*h*.

## 4- Studies on Exergy Analysis

Ozgener and Ozgener (2010) conducted an exergy analysis of an Earth-to-Air Heat Exchanger (EAHE) and applied it to a real-life system at the Solar Energy Institute of Izmir in Turkey. They independently calculated the exergetic efficiencies for each component. Their study reported that the maximum daily thermal coefficient of performance (COP) for the system was 6.18, with an average COP of 4.74 [23]. In the same year, they performed another study applying the EAHE system for greenhouse cooling under different operational conditions. They conducted an exergoeconomic analysis and found that the total exergy destruction for the system ranged from 0.26 kW to 2.5 kW. The maximum reported COP was 11.96, with an average COP of 5.89. Additionally, the overall exergetic efficiency, based on the product-to-fuel ratio, was 56.9% [24]. The results for greenhouse cooling were further explored in another study by Ozgener and Ozgener in 2010. When the minimum environmental conditions were 13.1°C temperature and 32% relative humidity, the maximum calculated COP was 6.42, with an average of 5.16 [25]. In 2011, Ozgener and Ozgener proposed an optimized design for a closed-loop EAHE system using exergoeconomic analysis [26]. In the same year, Yildiz et al. (2011) conducted an experimental study on the exergetic performance of an EAHE system enhanced with a photovoltaic solar system for greenhouse cooling. Their findings indicated that the system successfully cooled the target space. With an outside temperature of 15°C, the average exergetic efficiency was calculated as 56.3% for the heat exchanger and 4.94% for the photovoltaic system. Additionally, the system achieved an energy savings of 2.84 kW, equivalent to 31%, [27].

## 5. Conclusion

The Earth-to-Air Heat Exchanger (EAHE) system emerges as a promising solution for sustainable cooling and heating applications, particularly in residential and agricultural settings such as greenhouses. Its ability to integrate seamlessly with other green technologies and systems further amplifies its appeal. The system's low energy requirements, minimal environmental impact, and compatibility with renewable energy sources make it an exemplary choice for reducing reliance on conventional, energy-intensive methods. Exergoeconomic analyses provide robust evidence of the system's cost-effectiveness and energy efficiency, making it an economically viable solution for diverse climates and operational contexts. This study aimed to consolidate the growing body of research on EAHE systems and highlight their potential as a cornerstone of passive cooling and sustainable energy strategies. The findings emphasize that EAHE systems not only contribute to energy conservation but also offer a scalable and environmentally friendly alternative for future energy demands.